\documentclass[reqno,a4paper]{amsart}
\usepackage{amssymb,setspace}
\usepackage{ifpdf}
\ifpdf
 \usepackage[hyperindex,pagebackref]{hyperref}%
\else
 \expandafter\ifx\csname dvipdfm\endcsname\relax
 \usepackage[hypertex,hyperindex,pagebackref]{hyperref}
 \else
 \usepackage[dvipdfm,hyperindex,pagebackref]{hyperref}
 \fi
\fi
\allowdisplaybreaks[4]
\numberwithin{equation}{section}
\theoremstyle{plain}
\newtheorem{thm}{Theorem}[section]
\newtheorem{lem}{Lemma}[section]
\theoremstyle{remark}
\newtheorem{rem}{Remark}[section]
\theoremstyle{definition}
\newtheorem{dfn}{Definition}[section]
\DeclareMathOperator{\td}{d\mspace{-1mu}}
\newcommand{\cmdeg}[1]{\sideset{}{_\mathrm{cm}^{#1}}\deg}

\begin{document}

\title[Bessel functions and completely monotonic degree]
{Properties of modified Bessel functions and completely monotonic degrees of differences between exponential and trigamma functions}

\author[F. Qi]{Feng Qi}
\address{Institute of Mathematics, Henan Polytechnic University, Jiaozuo City, Henan Province, 454010, China}
\email{\href{mailto: F. Qi <qifeng618@gmail.com>}{qifeng618@gmail.com}, \href{mailto: F. Qi <qifeng618@hotmail.com>}{qifeng618@hotmail.com}, \href{mailto: F. Qi <qifeng618@qq.com>}{qifeng618@qq.com}}
\urladdr{\url{http://qifeng618.wordpress.com}}



\begin{abstract}
In the paper, the author establishes inequalities, monotonicity, convexity, and unimodality for functions concerning the modified Bessel functions of the first kind and compute the completely monotonic degrees of differences between the exponential and trigamma functions.
\end{abstract}

\thanks{This paper was typeset using \AmS-\LaTeX}

\maketitle

\tableofcontents

\section{Introduction}

For smoothness and consecution, we split this section into six subsections.

\subsection{An inequality and complete monotonicity}
In~\cite[Lemma~2]{Yang-Fan-2008-Dec-simp.tex}, the inequality
\begin{equation}\label{e-1-t-1}
\psi'(t)<e^{1/t}-1
\end{equation}
on $(0,\infty)$ was obtained and applied, where $\psi(t)$ stands for the digamma function which may be defined by the logarithmic derivative
\begin{equation*}
\psi(t)=[\ln\Gamma(t)]'=\frac{\Gamma'(t)}{\Gamma(t)}
\end{equation*}
and $\Gamma(t)$ is the classical Euler gamma function which may be defined for $\Re z>0$ by
\begin{equation*}
\Gamma(z)=\int^\infty_0t^{z-1} e^{-t}\td t.
\end{equation*}
The derivatives $\psi'(z)$ and $\psi''(z)$ are respectively called the tri- and tetra-gamma functions.
\par
In~\cite[Theorem~3.1]{QiBerg.tex}, \cite[Theorem~1.1]{simp-exp-degree-revised.tex}, and~\cite{simp-exp-degree-new.tex}, among other things, the inequality~\eqref{e-1-t-1} was generalized to a complete monotonicity respectively by different and elementary approaches, which reads that the difference
\begin{equation}\label{alpha-exp=psi-eq}
h(t)=e^{1/t}-\psi'(t)
\end{equation}
is completely monotonic, that is, $(-1)^{k-1}h^{(k-1)}(t)\ge0$ for $k\in\mathbb{N}$, on $(0,\infty)$ and
\begin{equation}\label{h(t)-limit=1}
\lim_{t\to\infty}h(t)=1.
\end{equation}

\subsection{Integral representations}

In~\cite[Theorem~1.2]{simp-exp-degree-revised.tex}, among other things, it was established that
\begin{equation}\label{exp=k=degree=k+1-int-bes}
e^{1/z}-\sum_{m=0}^k\frac{1}{m!}\frac1{z^m}=\frac1{z^{k+1}}\biggl[\frac1{(k+1)!}+\int_0^\infty \frac{I_{k+2} \bigl(2\sqrt{t}\,\bigr)}{t^{(k+2)/2}} e^{-zt}\td t\biggr]
\end{equation}
and
\begin{equation}\label{exp=k=degree=k+1-int}
e^{1/z}-\sum_{m=0}^k\frac{1}{m!}\frac1{z^m}=\frac1{k!(k+1)!}\int_0^\infty {}_1F_2(1;k+1,k+2;t)t^k e^{-zt}\td t
\end{equation}
for $\Re z>0$ and $k\in\{0\}\cup\mathbb{N}$, where the modified Bessel function of the first kind
\begin{equation}\label{I=nu(z)-eq}
I_\nu(z)= \sum_{k=0}^\infty\frac1{k!\Gamma(\nu+k+1)}\biggl(\frac{z}2\biggr)^{2k+\nu}
\end{equation}
for $\nu\in\mathbb{R}$ and $z\in\mathbb{C}$, the hypergeometric series
\begin{equation}
{}_pF_q(a_1,\dotsc,a_p;b_1,\dotsc,b_q;x)=\sum_{n=0}^\infty\frac{(a_1)_n\dotsm(a_p)_n} {(b_1)_n\dotsm(b_q)_n}\frac{x^n}{n!}
\end{equation}
for $b_i\notin\{0,-1,-2,\dotsc\}$, and the shifted factorial $(a)_0=1$ and
\begin{equation}
(a)_n=a(a+1)\dotsm(a+n-1)
\end{equation}
for $n>0$ and any real or complex number $a$. This gives an answer to an open problem posed in~\cite{exp-reciprocal-cm-IJOPCM.tex}. See also~\cite[Chapter~6]{Zhang-Xiao-Jing-Thesis.tex}.
\par
By virtue of~\eqref{exp=k=degree=k+1-int-bes} or~\eqref{exp=k=degree=k+1-int} for $k=0$ and the well known formula
\begin{equation}\label{polygamma}
\psi^{(n)}(z)=(-1)^{n+1}\int_0^{\infty}\frac{u^n}{1-e^{-u}}e^{-zu}\td u
\end{equation}
for $n=1$, see~\cite[p.~260, 6.4.1]{abram}, we can easily derive that
\begin{equation}\label{h(t)-int-repres}
h(z) =1+\int_0^{\infty}\biggl[\frac{I_1\bigl(2\sqrt{u}\,\bigr)}{\sqrt{u}\,} -\frac{u}{1-e^{-u}}\biggr]e^{-zu}\td u
\end{equation}
for $\Re z>0$. See the equation~(4.3) in~\cite{simp-exp-degree-revised.tex}.

\subsection{Lower bounds for a modified Bessel function of the first kind}

By the complete monotonicity obtained in~\cite[Theorem~3.1]{QiBerg.tex} and~\cite[Theorem~1.1]{simp-exp-degree-revised.tex} for $h(t)$, by the integral representation~\eqref{h(t)-int-repres}, and by Lemma~\ref{Bernstein-Widder-Theorem-12b} below, it was deduced in~\cite[Theorem~7.1]{QiBerg.tex} that
\begin{equation}\label{I=1-exp-ineq}
\alpha I_1(x)>\frac{(x/2)^3}{1-e^{-(x/2)^2}}
\quad\text{and}\quad
I_1(x)\ge \frac{\frac1\beta\bigl(\frac{x}2\bigr)^3} {1-\exp\bigl[-\frac1\beta\bigl(\frac{x}2\bigr)^2\bigr]}
\end{equation}
on $(0,\infty)$ if and only if $\alpha\ge1$ and $\beta\ge1$. More strongly, it was discovered in~\cite[Theorem~5.1]{QiBerg.tex} that, when $\beta\ge1$, the function
\begin{equation}\label{F-beta-(u)-dfn}
F_{\beta}(u)=\frac{u}{1-e^{-u}}\frac{\sqrt{\beta u}\,}{I_1\bigl(2\sqrt{\beta u}\,\bigr)}
\end{equation}
is decreasing on $(0,\infty)$; when $0<\beta<1$, it is unimodal (that is, it has a unique maximum) and $\frac1{F_{\beta}(u)}$ is convex on $(0,\infty)$.

\subsection{Necessary and sufficient conditions}

For $\alpha,\beta>0$, let
\begin{equation}\label{beta-psi-exp}
h_{\alpha,\beta}(t)=\alpha e^{\beta/t}-\psi'(t)
\end{equation}
on $(0,\infty)$. In~\cite{QiBerg.tex}, among other things, the following necessary and sufficient conditions for the function $h_{\alpha,\beta}(t)$ to be completely monotonic on $(0,\infty)$ were obtained.

\begin{thm}[{{\cite[Theorem~4.1]{QiBerg.tex}}}]\label{alpha-beta-psi-exp-thm-orig}
The function $h_{1,\beta}(t)$ is completely monotonic on $(0,\infty)$ if and only if $\beta\ge1$.
\par
If $\beta\ge1$ and $\alpha\beta\ge1$, the function $h_{\alpha,\beta}(t)$ is completely monotonic on $(0,\infty)$.
\par
A necessary condition for the function $h_{\alpha,\beta}(t)$ to be completely monotonic on $(0,\infty)$ is $\alpha\beta\ge1$.
\par
If $0<\beta<1$, the condition
\begin{equation}\label{beta<1-condition}
\alpha\beta\ge\max_{u\in(0,\infty)}F_{\beta}(u)>1
\end{equation}
is necessary and sufficient for $h_{\alpha,\beta}(t)$ to be completely monotonic on $(0,\infty)$, where
\begin{equation}\label{F-beta-2-limits}
\lim_{u\to0^+}F_{\beta}(u)=1\quad \text{and}\quad \lim_{u\to\infty}F_{\beta}(u)=0
\end{equation}
for all $\beta>0$.
\end{thm}

\subsection{Completely monotonic degree}

The notion ``completely monotonic degree'' was created in~\cite[Definition~1]{psi-proper-fraction-degree-two.tex}, which may be regarded as a slight but essential modification of~\cite[Definition~1.5]{Koumandos-Pedersen-09-JMAA}.
This definition may be further modified as follows.

\begin{dfn}\label{x-degree-dfn}
Let $f(x)$ be a completely monotonic function on $(0,\infty)$ and denote $f(\infty)=\lim_{x\to\infty}f(x)$. If for some $r\in\mathbb{R}$ the function $x^r[f(x)-f(\infty)]$ is completely monotonic on $(0,\infty)$ but $x^{r+\varepsilon}[f(x)-f(\infty)]$ is not for any positive number $\varepsilon>0$, then we say that the number $r$ is the completely monotonic degree of $f(x)$ with respect to $x\in(0,\infty)$; if for all $r\in\mathbb{R}$ each and every $x^r[f(x)-f(\infty)]$ is completely monotonic on $(0,\infty)$, then we say that the completely monotonic degree of $f(x)$ with respect to $x\in(0,\infty)$ is $\infty$.
\end{dfn}

In~\cite[p.~9890]{psi-proper-fraction-degree-two.tex}, the notation $\cmdeg{x}[f(x)]$ was designed to denote the completely monotonic degree $r$ of $f(x)$ with respect to $x\in(0,\infty)$. We can redevelop the above Definition~\ref{x-degree-dfn} as follows. If $f:(0,\infty)\to [0,\infty)$ is a $C^\infty$-function, then
\begin{equation}
\cmdeg{x}[f(x)]=\sup\{r\in\mathbb{R} \mid \text{$x^r[f(x)-f(\infty)]$ is completely monotonic}\}.
\end{equation}
It is clear that the completely monotonic degree $\cmdeg{x}[f(x)]$ of any completely monotonic function $f(x)$ on $(0,\infty)$ is at leat $0$.
\par
We claim that the completely monotonic degree $\cmdeg{x}[f(x)]$ equals $\infty$ if and only if $f(x)$ is nonnegative and identically constant.
It is clear that $\cmdeg{x}[0]=\infty$, as defined in~\cite[p.~9891, (4)]{psi-proper-fraction-degree-two.tex}.
Conversely, if $\cmdeg{x}[f(x)]=\infty$, then $x^r[f(x)-f(\infty)]$ is always completely monotonic on $(0,\infty)$ for any number $r\ge0$. This means that
\begin{align*}
\{x^r[f(x)-f(\infty)]\}'&=rx^{r-1}[f(x)-f(\infty)]+x^rf'(x)\\
&=x^{r-1}\{r[f(x)-f(\infty)]+xf'(x)\}\\
&\le0,
\end{align*}
that is,
\begin{equation}\label{rf(x)+xf(x)}
r[f(x)-f(\infty)]+xf'(x)\le0
\end{equation}
is valid on $(0,\infty)$ for all $r\ge0$. A result on~\cite[p.~98]{Dubourdieu} asserts that for a completely monotonic function $f$ on $(0,\infty)$ the strict inequality $(-1)^{k-1}f^{(k-1)}(t)>0$ for $k\in\mathbb{N}$ holds unless $f(x)$ is constant. This implies that, if $f(x)$ is not identically constant, then $f(x)-f(\infty)>0$ and $f'(x)<0$ on $(0,\infty)$. Consequently, the inequality~\eqref{rf(x)+xf(x)} may be rearranged as
$$
r<-\frac{xf'(x)}{f(x)-f(\infty)},\quad x\in(0,\infty).
$$
This leads to a contradiction to the arbitrariness of $r\ge0$. As a result, it holds that the function $f(x)$ is identically constant on $(0,\infty)$.

\subsection{Main results of this paper}

Since the complete monotonicity of $h(t)$ and the limit~\eqref{h(t)-limit=1} have been verified in~\cite[Theorem~3.1]{QiBerg.tex} and~\cite[Theorem~1.1]{simp-exp-degree-revised.tex}, we naturally consider to compute the completely monotonic degrees of the completely monotonic function
\begin{equation}\label{H-alpha-beta(t)-dfn}
H_{\alpha,\beta}(t)=h_{\alpha,\beta}(t)-\alpha
\end{equation}
on $(0,\infty)$.
\par
Our main results in this paper are the following theorems in sequence.

\begin{thm}\label{lower-beseel=5-thm}
When $1\le k\le 5$, the inequality
\begin{equation}\label{conj-bessel-ineq}
\frac{I_k\bigl(2\sqrt{u}\,\bigr)}{u^{k/2}}\ge\biggl(\frac{u}{1-e^{-u}}\biggr)^{(k-1)}
\end{equation}
is valid on $(0,\infty)$.
\par
When $\beta\ge1$, the function
\begin{equation}\label{G-beta-(u)-dfn}
G_{\beta}(u)=\frac{\beta u}{I_2\bigl(2\sqrt{\beta u}\,\bigr)}\biggl(\frac{u}{1-e^{-u}}\biggr)'
\end{equation}
is decreasing on $(0,\infty)$; when $0<\beta<1$, it is unimodal and $\frac1{G_{\beta}(u)}$ is convex on $(0,\infty)$.
\end{thm}

\begin{thm}\label{H(t)-degree=4-thm}
Let $\alpha,\beta>0$.
\begin{enumerate}
\item
If $(\alpha,\beta)=(1,1)$, then
\begin{equation}\label{H(t)degree=4}
\cmdeg{t}[H_{1,1}(t)]=4;
\end{equation}
\item
if $\beta>1$, then
\begin{equation}
\cmdeg{t}\bigl[H_{1/\beta,\beta}(t)\bigr]=2;
\end{equation}
\item
if $\alpha\beta>1$ and $\beta\ge1$, or if $\alpha\beta>1$, $0<\beta<1$, and
$$
\alpha\beta^2\ge\max_{u\in(0,\infty)}\{G_\beta(u)\},
$$
then we have
\begin{equation}\label{H-alpha-beta(t)degree=1}
\cmdeg{t}[H_{\alpha,\beta}(t)]=1.
\end{equation}
\end{enumerate}
\end{thm}
\section{Lemmas}

We need the following lemmas.

\begin{lem}[{{\cite[Theorem~1 and Lemma~1]{Constantinescu-Gen-Math-04}}}] \label{Const-Gen-Math-04-ineq-exp}
Let
\begin{equation}
c_k(n)=\binom{n}{k}\frac{(2n-k)!}{n!}\quad \text{and}\quad Q_n(x)=\sum_{k=0}^nc_k(n)x^{n-k}.
\end{equation}
If $m,n\in\mathbb{N}$, then the polynomial $Q_{2n}$ has no real root, the polynomial $Q_{2n+1}$ has a unique real root in $(-\infty,0)$, and
\begin{equation}
\frac{Q_{2n}(x)}{Q_{2n}(-x)}<e^{1/x}<-\frac{Q_{2m+1}(x)}{Q_{2m+1}(-x)},\quad x\ge\frac1{2(m+1)}.
\end{equation}
\end{lem}

\begin{lem}[{{\cite{Cargo-Shisha-1966} and~\cite[p.~227, 3.3.27]{mit}}}]\label{mit-3.3.27-ineq}
Let
$
P(x)=\sum_{k=0}^na_kx^k
$
be a real polynomial of degree $n\ge0$. Then the inequality
\begin{equation}
\min\{b_k\mid0\le k\le n\}\le P(x)\le\max\{b_k\mid0\le k\le n\}
\end{equation}
holds for $0\le x\le1$, where
$
b_k=\sum_{\ell=0}^ka_\ell\frac{\binom{k}{\ell}}{\binom{n}{\ell}}
$
for $0\le k\le n$.
\end{lem}

\begin{lem}[\cite{Biernacki-Krzyz-1955}]\label{Ponnusamy-Vuorinen-97-lem}
Let $a_k$ and $b_k$ for $k\in\{0\}\cup\mathbb{N}$ be real numbers and the power series
\begin{equation}
A(x)=\sum_{k=0}^\infty a_kx^k\quad\text{and}\quad B(x)=\sum_{k=0}^\infty
b_kx^k
\end{equation}
be convergent on $(-R,R)$ for some $R>0$. If $b_k>0$ and the ratio $\frac{a_k}{b_k}$ is \textup{(}strictly\textup{)} increasing for $k\in\mathbb{N}$, then the function $\frac{A(x)}{B(x)}$ is also \textup{(}strictly\textup{)} increasing on $(0,R)$.
\end{lem}

\begin{lem}[{{\cite[p.~377, 9.7.1]{abram}}}]\label{BesselI-Asympt-lem}
When $\nu$ is fixed, $|z|$ is large, and $\mu=4\nu^2$,
\begin{equation}
I_\nu(z)\sim\frac{e^z}{\sqrt{2\pi z}} \Biggl\{1+\sum_{\ell=1}^\infty \frac{(-1)^\ell}{\ell!(8z)^\ell} \prod_{j=1}^\ell\bigl[\mu-(2j-1)^2\bigr]\Biggr\},\quad |\arg z|<\frac\pi2.
\end{equation}
\end{lem}

\begin{lem}[{{\cite[p.~260, 6.4.11]{abram}}}]
For $|\arg z|<\pi$,
\begin{equation}\label{asymptotic-polypsi}
\psi^{(n)}(z)=(-1)^{n-1}\biggl[\frac{(n-1)!}{z^n}+\frac{n!}{2z^{n+1}}+\sum_{k=1}^\infty B_{2k}\frac{(2k+n-1)!}{(2k)!z^{2k+n}}\biggr],
\end{equation}
where $B_n$ for $n\ge0$ stand for Bernoulli numbers which may be generated by
\begin{equation}\label{Bernoulli-numbers-dfn}
\frac{x}{e^x-1}=\sum_{n=0}^\infty \frac{B_n}{n!}x^n =1-\frac{x}2+\sum_{j=1}^\infty B_{2j}\frac{x^{2j}}{(2j)!}, \quad\vert x\vert <2\pi.
\end{equation}
\end{lem}

\begin{lem}[{{\cite[p.~161, Theorem~12b]{widder}}}]\label{Bernstein-Widder-Theorem-12b}
A necessary and sufficient condition for $f(x)$ to be completely monotonic on $(0,\infty)$ is that
\begin{equation} \label{berstein-1}
f(x)=\int_0^\infty e^{-xt}\td\mu(t),
\end{equation}
where $\mu$ is a positive measure on $[0,\infty)$ such that the integral converges on $(0,\infty)$.
\end{lem}

\section{Proof of Theorem~\ref{lower-beseel=5-thm}}

For proving Theorem~\ref{H(t)-degree=4-thm}, we need at first to verify Theorem~\ref{lower-beseel=5-thm} as follows.

\begin{proof}[Proof of the inequality~\eqref{conj-bessel-ineq}]
Taking $\nu=5$ and $z=2\sqrt{u}\,$ in~\eqref{I=nu(z)-eq} lead to
\begin{equation}
\frac{I_5\bigl(2\sqrt{u}\,\bigr)}{u^{5/2}}=\sum_{k=0}^\infty\frac{u^k}{k!(k+5)!} \ge\sum_{k=0}^1\frac{u^k}{k!(k+5)!} =\frac{6+u}{720}.
\end{equation}
Hence, in order to prove the inequality~\eqref{conj-bessel-ineq} for $k=5$, it is sufficient to show
\begin{equation*}
\frac{u+6}{720}\ge\biggl(\frac{u}{1-e^{-u}}\biggr)^{(4)}
=\frac{e^u \bigl[(u-4)e^{3u}+(11u-12)e^{2u}+(11u+12)e^u+u+4\bigr]}{(e^u-1)^5}
\end{equation*}
which can be rewritten as
\begin{equation*}
\begin{split}
F_1(u)&\triangleq(u+6)(e^u-1)^5-720e^u \bigl[(u-4)e^{3u}+(11u-12)e^{2u}+(11u+12)e^u+u+4\bigr]\\
&\ge0.
\end{split}
\end{equation*}
\par
By direct calculations, we have
\begin{align*}
F_1'(u)&=(31+5 u)e^{5 u}+25(427-116 u) e^{4 u}+(18190-23730 u)e^{3 u}\\
&\quad-10 (2533+1586 u) e^{2 u} -5 (713+143 u) e^u-1,\\
F_1''(u)&=5 e^u \bigl[(32+5 u)e^{4 u}+40 (199-58 u) e^{3 u}+(6168-14238 u)e^{2 u}\\
&\quad-8 (1663+793 u) e^u-143 u-856\bigr]\\
&\triangleq 5 e^u F_2(t),\\
F_2'(u)&=(133+20 u)e^{4 u} +40 (539-174 u) e^{3 u}-6 (317+4746 u) e^{2 u}\\
&\quad-8 (2456+793 u) e^u-143,\\
F_2''(u)&=8 e^u \bigl[ (69+10 u)e^{3 u}+(7215-2610 u)e^{2 u} -3 (1345+2373 u) e^u\\
&\quad-793 u-3249\bigr]\\
&\triangleq8 e^uF_3(u),
\end{align*}
and
\begin{equation*}
F_3(0)=F_2''(0)=F_2'(0)= F_2(0)=F_1''(0)=F_1'(0)=F_1(0)=0.
\end{equation*}
As a result, when $F_3(u)\ge0$ on $(0,\infty)$, it follows that $F_1(u)\ge0$ on $(0,\infty)$.
\par
It is easy to see that the function $F_3(u)$ can be rearranged as
\begin{equation}
\begin{split}\label{F=3(u)-dfn}
F_3(u)&=10 u\bigl(e^{u}-261\bigr) e^{2 u}+3 \bigl(23 e^{2 u}-2373 u-1345\bigr) e^u\\
&\quad+7215 e^{2 u}-793 u-3249
\end{split}
\end{equation}
and that
\begin{enumerate}
\item
the term $e^u-261$ is positive when $u>\ln261=5.56\dotsc$,
\item
the term $23 e^{2 u}-2373 u-1345$ has a unique minimum at $u=\frac12\ln\frac{2373}{46}=1.97\dotsc$ on $(0,\infty)$ and equals $23\bigl(e^6-368\bigr)=814.86\dotsc$ at the point $u=3$,
\item
and $7215 e^{2u}-(3249+793u)$ has a unique minimum at $u=-\frac12\ln\frac{1110}{61}$ on $(-\infty,\infty)$ and is positive on $(0,\infty)$.
\end{enumerate}
Consequently, when $u\ge6$, the function $F_3(u)$ is positive.
\par
Applying Lemma~\ref{Const-Gen-Math-04-ineq-exp} to $m=2$ and $n=3$ derives
\begin{multline*}
\frac{x^6+42x^5+840x^4+10080x^3+75600x^2+332640x+665280}{x^6-42x^5+840x^4-10080x^3+75600 x^2-332640x+665280}<e^x\\
<\frac{x^5+30x^4+420x^3+3360x^2+15120x+30240}{30240-15120x+3360x^2-420x^3+30x^4-x^5}, \quad0<x\le6
\end{multline*}
and the function $F_3(u)$ can be written as
\begin{align*}
F_3(u)&=69e^{3u}+7215e^{2u}-4035e^u-3249+\bigl(10e^{3u}-2610e^{2u}-7119e^u-793\bigr)u\\
&>69\biggl(\frac{665280+332640 u+75600 u^2+10080 u^3+840 u^4+42 u^5+u^6}{665280-332640 u+75600 u^2-10080 u^3+840 u^4-42 u^5+u^6}\biggr)^3\\
&\quad+7215\biggl(\frac{665280+332640 u+75600 u^2+10080 u^3+840 u^4+42 u^5+u^6}{665280-332640 u+75600 u^2-10080 u^3+840 u^4-42 u^5+u^6}\biggr)^2\\
&\quad-4035\biggl(\frac{30240+15120 u+3360 u^2+420 u^3+30 u^4+u^5}{30240-15120 u+3360 u^2-420 u^3+30 u^4-u^5}\biggr)-3249\\
&\quad+u\biggl[10\biggl(\frac{665280+332640 u+75600 u^2+10080 u^3+840 u^4+42 u^5+u^6}{665280-332640 u+75600 u^2-10080 u^3+840 u^4-42 u^5+u^6}\biggr)^3\\
&\quad-2610\biggl(\frac{30240+15120 u+3360 u^2+420 u^3+30 u^4+u^5}{30240-15120 u+3360 u^2-420 u^3+30 u^4-u^5}\biggr)^2\\
&\quad-7119\biggl(\frac{30240+15120 u+3360 u^2+420 u^3+30 u^4+u^5}{30240-15120 u+3360 u^2-420 u^3+30 u^4-u^5}\biggr)-793\biggr]\\
&\triangleq\frac{6uF_4(u)}{F_5(u)}
\end{align*}
on $(0,6]$, where
\begin{align*}
F_4(u)&=621 u^{28}-94751 u^{27}+6068010 u^{26}-218010408 u^{25}+4603805304u^{24}\\
&\quad-34479103680 u^{23}-1412513172000 u^{22}+69509082484800u^{21}\\
&\quad-1835940264439680 u^{20}+36018288433370880u^{19}\\
&\quad-572728070517926400u^{18}+7670777548637952000u^{17}\\
&\quad-88289462254568601600u^{16}+882241752079928217600u^{15}\\
&\quad-7678232793596974694400 u^{14}+58018545121380802560000u^{13}\\
&\quad-376773142967398969344000 u^{12}+2061377592729654140928000u^{11}\\
&\quad-9156137875572402634752000u^{10}+30521365267364424843264000 u^9\\
&\quad-59514097618800165519360000u^8-48948451585366441328640000 u^7\\
&\quad+852510196971380523663360000u^6-3082810530742053482004480000u^5\\
&\quad+5063190015183760203448320000 u^4-350858087962497987379200000u^3\\
&\quad-10582859071799067200716800000u^2\\
&\quad+8481790289232945532108800000 u+4038947756777593110528000000
\end{align*}
and
\begin{align*}
F_5(u)&=\bigl(u^5-30 u^4+420 u^3-3360 u^2+15120 u-30240\bigr)^2 \\
&\quad\times\bigl(u^6-42 u^5+840 u^4-10080 u^3+75600 u^2-332640 u+665280\bigr)^3.
\end{align*}
\par
By Lemma~\ref{Const-Gen-Math-04-ineq-exp}, it follows that the function
\begin{equation*}
F_6(u)=u^5Q_5\biggl(-\frac1u\biggr)=u^5-30 u^4+420 u^3-3360 u^2+15120 u-30240
\end{equation*}
has a unique positive zero. Since $F_6(6)=-864$ and $F_6(8)=608$, the function $F_6(u)$ is negative on $(0,6)$. By Lemma~\ref{Const-Gen-Math-04-ineq-exp} once again, the function
\begin{equation*}
u^6Q_6\biggl(-\frac1u\biggr)=u^6-42 u^5+840 u^4-10080 u^3+75600 u^2-332640 u+665280
\end{equation*}
has no any zero. In a word, the function $F_5(u)>0$ on $(0,6)$.
\par
A direct computation shows that the sequence $\{b_k\mid0\le k\le n\}$ in Lemma~\ref{mit-3.3.27-ineq} applied to $F_4(u)$ are
\begin{align*}
b_0&=4038947756777593110528000000, \\
b_1&=4341868838535912593817600000, \\
b_2&=4616792938622805973401600000, \\
b_3&=\frac{63226968447497701058150400000}{13}, \\
b_4&=\frac{66072098066989847041120665600}{13}, \\
b_5&=\frac{68558326435040659929169920000}{13}, \\
b_6&=\frac{1625932364702092812073304064000}{299}, \\
b_7&=\frac{18338550757492925804643090432000}{3289}, \\
b_8&=\frac{18707897201998227336080916480000}{3289}, \\ b_9&=\frac{18996421971796176266488102256640}{3289},\\
b_{10}&=\frac{364943527309321738804310099558400}{62491}, \\ b_{11}&=\frac{33414027455107862233493347123200}{5681}, \\
b_{12}&=\frac{570021720890042639974424894668800}{96577}, \\ b_{13}&=\frac{43851581059726626552866155622400}{7429}, \\
b_{14}&=\frac{43715751783929803366617869881344}{7429}, \\ b_{15}&=\frac{43449666570024250746295234195200}{7429}, \\
b_{16}&=\frac{29463945191749248714150342727680}{5083}, \\ b_{17}&=\frac{32549053609077425759197656000000}{5681}, \\
b_{18}&=\frac{352985157442262967317198674456320}{62491}, \\ b_{19}&=\frac{91354286436246580611619508370624}{16445}, \\
b_{20}&=\frac{17926195304232242164107747202432}{3289}, \\ b_{21}&=\frac{17548009393294868476730762864424}{3289}, \\
b_{22}&=\frac{67747076984066766537183527176}{13}, \\
b_{23}&=\frac{66030748681802854683680196472}{13}, \\
b_{24}&=\frac{4816825337954730130116871131704}{975}, \\
b_{25}&=\frac{62340864945822949546899440674}{13},\\
b_{26}&=\frac{292672542741083383435627679675}{63}, \\
b_{27}&=\frac{125766913766925341535184862941}{28}, \\
b_{28}&=4334548991696365872138512296.
\end{align*}
This means, by virtue of Lemma~\ref{mit-3.3.27-ineq}, that $F_4(u)$ is positive on $[0,1]$.
\par
The positivity of $F_4(u)$ on the interval $[1,2]$, $[2,3]$, $[3,4]$, $[4,5]$, or $[5,6]$ can be respectively transformed into the positivity of the function
\begin{gather*}
F_5(u)=F_4(u+1), \quad F_6(u)=F_5(u+1)=F_4(u+2), \\
F_7(u)=F_6(u+1)=F_5(u+2)=F_4(u+3), \\
F_8(u)=F_7(u+1)=F_6(u+2)=F_5(u+3)=F_4(u+4),
\end{gather*}
or
\begin{equation*}
F_9(u)=F_8(u+1)=F_7(u+2)=F_6(u+3)=F_5(u+4)=F_4(u+5)
\end{equation*}
on the unit interval $[0,1]$, which can be respectively verified by Lemma~\ref{mit-3.3.27-ineq} as done in the proof of the positivity of the function $F_4(u)$ on $[0,1]$.
\par
In conclusion, the function $F_3(u)$, and so $F_1(u)$, is positive on $(0,\infty)$. This means that the inequality~\eqref{conj-bessel-ineq} for $k=5$ is valid on $(0,\infty)$.
\par
The inequality~\eqref{conj-bessel-ineq} for $1\le k\le4$ may be verified by similar arguments as above.
\end{proof}

\begin{proof}[Proof of monotonicity of the function~\eqref{G-beta-(u)-dfn}]
For simplicity, we consider
\begin{equation*}
\frac1{G_{\beta}(u)}=\frac{(e^u-1)^2}{(e^u-1-u)e^u} \frac{I_2\bigl(2\sqrt{\beta u}\,\bigr)}{\beta u}
\triangleq \frac{Q_{\beta}(u)}{P(u)},
\end{equation*}
where
\begin{equation*}
P(u)=(e^u-1-u)e^u=\sum_{k=2}^\infty\frac{2^k-k-1}{k!}u^k=\sum_{k=2}^\infty p_{k-2}u^k
\end{equation*}
and
\begin{equation*}
Q_{\beta}(u)=(e^u-1)^2 \frac{I_2\bigl(2\sqrt{\beta u}\,\bigr)}{\beta u}
=\sum_{k=2}^\infty\frac{2^k-2}{k!}u^k\sum_{k=0}^\infty\frac{\beta^k}{k!(k+2)!}u^k
=\sum_{k=2}^\infty q_{k-2}(\beta)u^k
\end{equation*}
with
\begin{equation*}
p_k=\frac{2^{k+2}-k-3}{(k+2)!}
\quad\text{and}\quad
q_k(\beta)=\frac{1}{(k+2)!} \sum_{\ell=0}^k\binom{k+2}{\ell}\frac{2^{k-\ell+2}-2}{(\ell+2)!}\beta^\ell
\end{equation*}
for $k\ge0$. Hence,
\begin{equation*}
\frac{Q_{\beta}(u)}{P(u)} =\frac{\sum_{k=0}^\infty q_k(\beta)u^k} {\sum_{k=0}^\infty p_ku^k}.
\end{equation*}
\par
When $\beta\ge1$, let
\begin{equation*}
c_k(\beta)=\frac{q_k(\beta)}{p_k} =\frac1{2^{k+2}-k-3} \sum_{\ell=0}^k\binom{k+2}{\ell}\frac{2^{k-\ell+2}-2} {(\ell+2)!}\beta^\ell
\end{equation*}
for $k\in\{0\}\cup\mathbb{N}$.
It is clear that
\begin{equation*}
c_0(\beta)=1,\quad c_1(\beta)=\frac{3+\beta}4,\quad \text{and}\quad c_2(\beta)=\frac{14+8\beta+\beta^2}{22}
\end{equation*}
satisfy $c_0(\beta)\le c_1(\beta)<c_2(\beta)$ for $\beta\ge1$.
For $k\ge 2$,
\begin{equation*}
c_{k+1}(\beta)- c_k(\beta)=\sum_{\ell=0}^k(b_{k+1,\ell}-b_{k,\ell})\frac{\beta^\ell}{(\ell+2)!} +\frac{\beta^{k+1}}{(k+1)!(2^{k+3}-k-4)},
\end{equation*}
where
\begin{equation*}
b_{k,\ell}=\frac{2^{k-\ell+2}-2}{2^{k+2}-k-3}\binom{k+2}{\ell}.
\end{equation*}
An easy computation yields
\begin{align*}
b_{k+1,0}-b_{k,0}&=-\frac{2\bigl(1+2^{k+1} k\bigr)} {(2^{k+2}-k-3)(2^{k+3}-k-4)}\\
&<0,\\
b_{k+1,1}-b_{k,1}&=\frac{2 \bigl\{1+2^k\bigl[2^{k+3}-\bigl(k^2+2k+6\bigr)\bigr]\bigr\}} {(2^{k+2}-k-3)(2^{k+3}-k-4)}\\
&\ge\frac{2\bigl\{1+2^k\bigl[2^3\bigl(1+k+k(k-1)/2\bigr) -\bigl(k^2+2k+6\bigr)\bigr]\bigr\}}{(2^{k+2}-k-3)(2^{k+3}-k-4)} \\
&=\frac{2\bigl[1+2^k\bigl(3k^2+2k+2\bigr)\bigr]}{(2^{k+2}-k-3)(2^{k+3}-k-4)}\\
&>0,
\end{align*}
and
\begin{multline*}
\frac{b_{k+1,0}-b_{k,0}}{2!}+\beta\frac{b_{k+1,1}-b_{k,1}}{3!}
\ge\frac{b_{k+1,0}-b_{k,0}}{2!}+\frac{b_{k+1,1}-b_{k,1}}{3!}\\
\begin{aligned}
&=\frac{2^{k}\bigl[2^{k+3}-\bigl(k^2+8k+6\bigr)\bigr]-2}{3(2^{k+2}-k-3) (2^{k+3}-k-4)}\\
&>\frac{2^{k}\bigl\{2^3[1+k+k(k-1)/2]-\bigl(k^2+8k+6\bigr)\bigr\}-2} {3(2^{k+2}-k-3) (2^{k+3}-k-4)}\\
&=\frac{2^{k}(3k-4)k+2^{k+1}-2} {3(2^{k+2}-k-3) (2^{k+3}-k-4)}\\
&>0
\end{aligned}
\end{multline*}
for $k\ge2$. Further, when $k\ge\ell\ge2$, we have
\begin{align*}
\frac{b_{k+1,\ell}-b_{k,\ell}}{\binom{k+2}{\ell}}
&=\frac{(k+3)\bigl(2^{k-\ell+3}-2\bigr)}{(k-\ell+3)(2^{k+3}-k-4)} -\frac{2^{k-\ell+2}-2}{2^{k+2}-k-3}\\
&=\frac{2^{k-\ell+2}-2}{2^{k+3}-k-4} \biggl[\frac{(k+3)\bigl(2^{k-\ell+3}-2\bigr)}{(k-\ell+3)(2^{k-\ell+2}-2)} -\frac{2^{k+3}-k-4}{2^{k+2}-k-3}\biggr]\\
&=\frac{2^{k-\ell+2}-2}{2^{k+3}-k-4} \biggl[\frac{2(k+3)}{k-\ell+3} \biggl(1+\frac1{2^{k-\ell+2}-2}\biggr) -\frac{2^{k+3}-k-4}{2^{k+2}-k-3}\biggr]\\
&\ge\frac{2^{k-\ell+2}-2}{2^{k+3}-k-4} \biggl[\frac{2(k+3)}{k+1} \biggl(1+\frac1{2^{k}-2}\biggr) -\frac{2^{k+3}-k-4}{2^{k+2}-k-3}\biggr]\\
&=\frac{2^{k-\ell+2}-2}{2^{k+3}-k-4} \biggl(\frac{k+3}{k+1} \frac{2^{k+1}-2}{2^{k}-2} -\frac{2^{k+3}-k-4}{2^{k+2}-k-3}\biggr)\\
&=\frac{\bigl(2^{k-\ell+2}-2\bigr) \bigl\{2^k\bigl[2^{k+4}-\bigl(k^2-k+22\bigr)\bigr]+2(k+5)\bigr\}} {(k+1)(2^{k}-2)(2^{k+2}-k-3)(2^{k+3}-k-4)}\\
&>0,
\end{align*}
where in the last line we used the inequality
\begin{equation*}
2^{k+4}>2^4\biggl[1+k+\frac{k(k-1)}2\biggr] =\bigl(k^2-k+22\bigr) +\bigl(7k^2+9k-6\bigr)>k^2-k+22.
\end{equation*}
Consequently, when $\beta\ge1$ and $k\ge2$,
\begin{equation*}
c_{k+1}(\beta)- c_k(\beta)>\frac{b_{k+1,0}-b_{k,0}}{2!}+\frac{b_{k+1,1}-b_{k,1}}{3!}\beta +\sum_{\ell=2}^k\frac{b_{k+1,\ell}-b_{k,\ell}}{(\ell+2)!}\beta^\ell>0.
\end{equation*}
Therefore, when $\beta\ge1$, we have $c_{k+1}(\beta)- c_k(\beta)>0$. Equivalently speaking, when $\beta\ge 1$, the sequence $c_k(\beta)$ is increasing with respect to $k\ge0$. From this and Lemma~\ref{Ponnusamy-Vuorinen-97-lem}, it follows that, when $\beta\ge1$, the function $\frac{Q_{\beta}(u)}{P(u)}$ is increasing on $(0,\infty)$. As a result, when $\beta\ge1$, the function $G_{\beta}(u)$ is decreasing on $(0,\infty)$. The proof of monotonicity of the function~\eqref{G-beta-(u)-dfn} is complete.
\end{proof}

\begin{proof}[Proof of unimodality and convexity of the function~\eqref{G-beta-(u)-dfn}]
By~\eqref{I=nu(z)-eq}, it is straightforward to obtain
\begin{equation}\label{Beseel-der-recur}
\frac{\td}{\td z}\biggl[\frac{I_\nu(z)}{(z/2)^\nu}\biggr]
=\frac{I_{\nu+1}(z)}{(z/2)^\nu}.
\end{equation}
Making use of~\eqref{Beseel-der-recur} and differentiating lead to
\begin{align*}
\frac{\td}{\td u}\biggl[\frac1{G_{\beta}(u)}\biggr]
&=\frac{e^u-1}{e^{u}(1-e^u+u)^2}\biggl\{(e^u-1) (e^u-1-u) \frac{\td}{\td u}\biggl[\frac{I_2\bigl(2\sqrt{\beta u}\,\bigr)}{\beta u}\biggr] \\
&\quad-[e^u (u-2)+u+2] \frac{I_2\bigl(2\sqrt{\beta u}\,\bigr)}{\beta u}\biggr\}\\
&=\frac1{e^{u}(e^u-1-u)^2}\biggl\{\beta(e^u-1)^2 (e^u-1-u) \frac{I_3\bigl(2\sqrt{\beta u}\,\bigr)}{(\beta u)^{3/2}} \\
&\quad-[e^u (u-2)+u+2] (e^u-1)\frac{I_2\bigl(2\sqrt{\beta u}\,\bigr)}{\beta u}\biggr\}\\
&\triangleq\frac{R_\beta(u)}{S(u)},
\end{align*}
where
\begin{align*}
S(u)=(e^u-1-u)^2e^u
&=u^4\sum_{k=0}^\infty\frac{3^{k+4}-(k+6)2^{k+4}+k^2+9k+21}{(k+4)!}u^k\\
&\triangleq u^4\sum_{k=0}^\infty \lambda_ku^k,\\
(e^u-1)^2 (e^u-1-u)&=\sum_{k=4}^\infty\frac{3^k-(k+6)2^{k-1}+2k+3}{k!}u^k,\\
[e^u (u-2)+u+2](e^u-1)&=\sum_{k=4}^\infty\frac{(k-4)2^{k-1}+4}{k!}u^k,
\end{align*}
and
\begin{align*}
R_\beta(u)&=\beta(e^u-1)^2 (e^u-1-u) \frac{I_3\bigl(2\sqrt{\beta u}\,\bigr)}{(\beta u)^{3/2}}
-[e^u(u-2)+u+2](e^u-1) \frac{I_2\bigl(2\sqrt{\beta u}\,\bigr)}{\beta u}\\
&=\sum_{k=4}^\infty\frac{3^k-(k+6)2^{k-1}+2k+3}{k!}u^k \sum_{k=0}^\infty\frac{\beta^{k+1}}{k!(k+3)!}u^k\\
&\quad-\sum_{k=4}^\infty\frac{(k-4)2^{k-1}+4}{k!}u^k \sum_{k=0}^\infty\frac{\beta^k}{k!(k+2)!}u^k\\
&=\sum_{k=0}^\infty\sum_{\ell=0}^k\binom{k+4}{\ell} \frac{3^{k-\ell+4}-(k-\ell+10)2^{k-\ell+3}+2(k-\ell)+11}{(k+4)!(\ell+3)!}\beta^{\ell+1}u^{k+4}\\
&\quad-u^4\sum_{k=0}^\infty\frac{u^k}{(k+4)!}\sum_{\ell=0}^k\binom{k+4}{\ell} \frac{(k-\ell)2^{k-\ell+3}+4}{(\ell+2)!}\\
&=u^4\sum_{k=0}^\infty\frac{u^k}{(k+4)!}\sum_{\ell=0}^k\binom{k+4}{\ell} \frac{\beta^\ell}{(\ell+3)!}\bigl\{\bigl[3^{k-\ell+4} -(k-\ell+10)2^{k-\ell+3} \\
&\quad +2(k-\ell)+11\bigr]\beta-(\ell+3)\bigl[(k-\ell)2^{k-\ell+3}+4\bigr]\bigr\}\\
&\triangleq u^4\sum_{k=0}^\infty \xi_k(\beta)u^k.
\end{align*}
When $0<\beta<1$, let $C_k(\beta)=\frac{\xi_k(\beta)}{\lambda_k}$, that is,
\begin{align*}
C_k(\beta) &=\frac1{U_k}\biggl[\frac{k+4}{2\times k!}\beta^{k+1} -2\bigl(2^{k+1}k+1\bigr)
+\sum_{\ell=1}^k\frac{(k+4)!} {\ell!(\ell+2)!(k-\ell+5)!}V_k(\ell)\beta^\ell\biggr]\\
&\triangleq\sum_{\ell=0}^{k+1}\theta_{k,\ell}\beta^\ell
\end{align*}
for $k\ge0$, where
\begin{equation*}
U_k=3^{k+4}-(k+6)2^{k+4}+k^2+9k+21
\end{equation*}
and
\begin{align*}
V_k(\ell)&=3^{k-\ell+5}\ell +2^{k-\ell+3}\bigl(\ell^2-17\ell-5k-k^2\bigr) -2\ell^2+2k\ell+17\ell-4k-20\\
&=(k-m)3^{m+5}+\bigl[m^2+(17-2k)m-22k\bigr]2^{m+3} -2 m^2+(2k-17)m+13 k-20\\
&\triangleq W_k(m)
\end{align*}
with $0\le m=k-\ell<k$.
It is not difficult to obtain that
\begin{gather*}
C_0(\beta)=\frac{\beta-1}{3},\quad C_1(\beta)=\frac{\beta^2+4\beta-4}{20},\quad C_2(\beta)=\frac{3 \beta^3+35 \beta^2+80 \beta-68}{520},\\
C_3(\beta)=\frac{\beta^4+22 \beta^3+150 \beta^2+256 \beta-168}{1872},
\end{gather*}
and
\begin{equation*}
C_4(\beta)=\frac{\beta^5+35 \beta^4+448 \beta^3+2268 \beta^2+3220 \beta-1548}{24444}.
\end{equation*}
Therefore, the differences
\begin{gather*}
C_1(\beta)-C_0(\beta)=\frac{3 \beta^2+8(1-\beta)}{60},\quad C_2(\beta)-C_1(\beta)=\frac{3\bigl[\beta^3+3 \beta^2+4(3-2\beta)\bigr]}{520},\\
C_3(\beta)-C_2(\beta)=\frac{5 \beta^4+56 \beta^3+120 \beta^2+32(12-5\beta)}{9360},\\
C_4(\beta)-C_3(\beta)=\frac{52 \beta^5+1141 \beta^4+8358 \beta^3+16086 \beta^2+24(1399-266\beta)}{1271088}
\end{gather*}
are all positive for $0<\beta<1$.
\par
For $k\ge4$ and $0<\beta<1$, we have
\begin{equation*}
C_{k+1}(\beta)-C_k(\beta)=\sum_{\ell=0}^{k+1}(\theta_{k+1,\ell}-\theta_{k,\ell})\beta^\ell +\theta_{k+1,k+2}\beta^{k+2}.
\end{equation*}
Since
\begin{align*}
U_k&>2^{k+4}+(k+4)2^{k+3}+\binom{k+4}{2}2^{k+2}-(k+6)2^{k+4}+k^2+9k+21\\
&=\bigl(k^2+3k-12\bigr)2^{k+1}+k^2+9k+21\\
&>0
\end{align*}
for $k\ge4$, we easily obtain that $\theta_{k+1,k+2}>0$ and
\begin{equation}\label{C=k+1=C-k>0-eq}
C_{k+1}(\beta)-C_k(\beta)>\sum_{\ell=0}^{k+1}(\theta_{k+1,\ell}-\theta_{k,\ell})\beta^\ell \end{equation}\label{c-difference-ineq}
for $k\ge4$ and $0<\beta<1$.
\par
The inequality
\begin{equation}\label{theta-0-compute}
\theta_{k+1,0}\ge\theta_{k,0}
\end{equation}
may be rewritten as
\begin{equation*}
3^{k+4} \bigl[(k-2)2^k+1\bigr]+2^k\bigl(192\times2^{k}-k^3-9 k^2-37 k-106\bigr)+k+5\ge0,
\end{equation*}
which may be deduced from
\begin{align*}
192\times2^{k}-k^3-9 k^2-37 k-106&>192\sum_{\ell=0}^3\binom{k}{\ell}-k^3-9 k^2-37 k-106\\
&=(31 k-9) k^2+123 k+86\\
&>0.
\end{align*}
Thus, the inequality~\eqref{theta-0-compute} must be valid for $k\ge2$.
\par
The inequality
\begin{equation}\label{theta=k>k-eq}
\theta_{k+1,\ell}\ge\theta_{k,\ell}
\end{equation}
for $k\ge4$ and $k\ge\ell\ge1$ can be rearranged as
\begin{equation*}
\frac{U_{k+1}}{(k+5)U_{k}}\le\frac{V_{k+1}(\ell)}{(k-\ell+6)V_k(\ell)} =\frac{W_{k+1}(m+1)}{(m+6)W_k(m)},
\end{equation*}
where $0\le m=k-\ell<k$. Furthermore, for $k\ge4$ and $0\le m\le k-2$, the inequality
\begin{equation}\label{W-m-decrease}
\frac{W_{k+1}(m+1)}{(m+6)W_k(m)}\ge\frac{W_{k+1}(m+2)}{(m+7)W_k(m+1)}
\end{equation}
may be rearranged as
\begin{equation}\label{M-m(k)-dfn}
M_m(k)\triangleq\mathcal{A}(m)k^2+\mathcal{B}(m)k+\mathcal{C}(m)\ge0,
\end{equation}
where
\begin{align*}
\mathcal{A}(m)&=\bigl(4m^3+86m^2+442m+276\bigr)2^{m+3} +\bigl(m^2+25m+150\bigr)4^{m+5}-2\bigl(4m^2\\
&\quad+40m +87\bigr)3^{m+5} +9^{m+6}+\bigl(m^2+m-102\bigr)2^{m+4} 3^{m+5}+4m^2+64m+249\\
&>\bigl(4m^3+86m^2+442m+276\bigr)2^{m+3} +4m^2+64m+249\\
&\quad+\bigl(m^2+25m+150\bigr)\bigl[3^{m+5}+(m+5)3^{m+4}\bigr]  +9^{m+6}\\
&\quad-2\bigl(4m^2+40m+87\bigr)3^{m+5}+\bigl(m^2+m-102\bigr)2^{m+4} 3^{m+5}\\
&=\bigl(4m^3+86m^2+442m+276\bigr)2^{m+3}+9^{m+6} +4m^2+64m+249\\
&\quad+\bigl(678+110 m+9 m^2+m^3\bigr)3^{m+4} +\bigl(m^2+m-102\bigr)2^{m+4} 3^{m+5}\\
&=\bigl(4m^3+86m^2+442m+276\bigr)2^{m+3}+\bigl(678+110 m+9 m^2+m^3\bigr)3^{m+4} \\
&\quad +3^{m+5}\bigl[3^{m+7}+2^{m+4}\bigl(m^2+m-102\bigr)\bigr] +4m^2+64m+249\\
&>\bigl(4m^3+86m^2+442m+276\bigr)2^{m+3}+\bigl(678+110 m+9 m^2+m^3\bigr)3^{m+4} \\
&\quad +3^{m+5}\Biggl[\sum_{\ell=0}^3\binom{m+7}{\ell}2^{m-\ell+7} +2^{m+4}\bigl(m^2+m-102\bigr)\Biggr]+4m^2+64m+249\\
&=\bigl(4m^3+86m^2+442m+276\bigr)2^{m+3}+\bigl(678+110 m+9 m^2+m^3\bigr)3^{m+4} \\
&\quad +\bigl(66+215 m+30 m^2+m^3\bigr)2^{3+m}3^{m+4} +4m^2+64m+249\\
&>0,\\
\mathcal{B}(m)&=2\bigl(8m^3+92m^2+282m+207\bigr)3^{m+5}-(2m+1)9^{m+6}-\bigl(6m^4+145m^3\\
&\quad+839m^2+592m-1524\bigr)2^{m+3}- \bigl(m^3+31m^2+234m+108\bigr)4^{m+5}\\
&\quad-\bigl(m^3+5m^2-96m-12\bigr)2^{m+3} 3^{m+6} -\bigl(8m^3+148m^2+794m+1065\bigr),
\end{align*}
and
\begin{align*}
\mathcal{C}(m)&=\bigl(2m^5+57m^4+388m^3+585m^2+480m+2988\bigr)2^{m+3}+ \bigl(m^4+36m^3 \\
&\quad+323m^2+12m-756\bigr)4^{m+4}+4m^4+84m^3+569m^2+1401m+1152\\
&\quad+\bigl(m^4+10m^3-99m^2+24m+252\bigr)2^{m+3} 3^{m+5}+m(m+1)9^{m+6}\\
&\quad -2\bigl(4m^4+52m^3+207m^2+327m+288\bigr)3^{m+5}.
\end{align*}
It is clear that $M_m(k)$ may be regarded as a quadratic polynomial of $k$ and it has a unique possible minimum point $-\frac{\mathcal{B}(m)}{2\mathcal{A}(m)}$, which, due to $k\ge4$ and $0\le m\le k-2$, should satisfy $-\frac{\mathcal{B}(m)}{2\mathcal{A}(m)}\ge m+2$.
But, the fact is that $-\frac{\mathcal{B}(m)}{2\mathcal{A}(m)}<m+2$, that is,
\begin{align*}
2(m&+2)\mathcal{A}(m)+\mathcal{B}(m)=\bigl[3^{m+8}+\bigl(m^3-3m^2-112 m-780\bigr)2^{m+3}\bigr]3^{m+5}\\
&\quad+\bigl(m^3+23m^2+166m+492\bigr)4^{m+5}-2\bigl(4m^2+52m+141\bigr)3^{m+5}\\
&\quad+\bigl(2m^4+43m^3+389m^2+1728m+2628\bigr)2^{m+3}-4m^2-40m-69\\
&>\Biggl[3^{8}\sum_{\ell=0}^3\binom{m}{\ell}2^{m-\ell}+\bigl(m^3-3m^2-112 m-780\bigr)2^{m+3}\Biggr]3^{m+5}\\
&\quad+\bigl(m^3+23m^2+166m+492\bigr)3^{m+5}-2\bigl(4m^2+52m+141\bigr)3^{m+5}\\
&\quad+\bigl(2m^4+43m^3+389m^2+1728m+2628\bigr)-4m^2-40m-69\\
&=\bigl(5136+29404m+6177m^2+2315m^3\bigr)2^{m-4}3^{m+5}+\bigl(m^3+15m^2\\
&\quad+62m+210\bigr)3^{m+5}+2m^4+43m^3+385m^2+1688m+2559\\
&>0.
\end{align*}
This contradiction shows that, when $k\ge4$ and $k\ge m+2\ge2$, the quantity $M_m(k)$ can be regarded as a quadratic polynomial of $k$ and it has no any minimum. Combining this with the fact that $\mathcal{A}(m)>0$ concludes that the quadratic polynomial $M_m(k)$ of $k$ is increasing with respect to $k$. A direct computation reveals that
\begin{gather*}
M_0(k)=3360\bigl(54-137k+74k^2\bigr), \quad M_1(k)=1568\bigl(6480-7306k+1909k^2\bigr),\\
M_2(k)=336\bigl(750942-549881k+95837k^2\bigr)
\end{gather*}
are positive for $k\ge4$ and that for $m\ge3$ and $k\ge m+2$
\begin{align*}
M_m(k)&\ge M_m(m+2)=2 \bigl\{9+\bigl(m^4+20m^3+155m^2+508m+780\bigr)2^{2m+7} \\
&\quad+\bigl(2m^4+61m^3+703m^2+3692m+7140\bigr)2^{m+2}+3^{m+5}\bigl[3^{m+7}\\
&\quad-\bigl(492+172m +29m^2+m^3\bigr)2^{m+2}-2\bigl(111+32m+2m^2\bigr)\bigr]\bigr\}\\
&>2 \Biggl\{9+\bigl(m^4+20m^3+155m^2+508m+780\bigr)2^{2m+7} +\bigl(2m^4+61m^3\\
&\quad+703m^2+3692m+7140\bigr)2^{m+2} +3^{m+5}\Biggl[3^{7}\sum_{\ell=0}^3\binom{m}{\ell}2^{m-\ell}\\
&\quad-\bigl(492+172m +29m^2+m^3\bigr)2^{m+2}
-2\bigl(111+32m+2m^2\bigr)\Biggr]\Biggr\}\\
&=2 \biggl\{\frac{3^{m+5}}{16}\bigl[665m^32^m +\bigl(331\times 2^m-64\bigr)m^2+4\bigl(893\times 2^m-256\bigr)m\\
&\quad+48\bigl(73\times2^m-74\bigr)\bigr]+\bigl(m^4+20m^3+155m^2+508m+780\bigr)2^{2m+7} \\
&\quad+\bigl(2m^4+61m^3+703m^2+3692m+7140\bigr)2^{m+2}+9\biggr\}\\
&>0.
\end{align*}
Accordingly, the inequality~\eqref{M-m(k)-dfn}, and so the inequality~\eqref{W-m-decrease}, holds for all $0\le m\le k-2$ and $k\ge4$. This means that the sequence $\frac{W_{k+1}(m+1)}{(m+6)W_k(m)}$ is decreasing with respect to $m$, and so that the sequence $\frac{V_{k+1}(\ell)}{(k-\ell+6)V_k(\ell)}$ is increasing with respect to $\ell$. Therefore, in order to show the inequality~\eqref{theta=k>k-eq} for $k\ge4$ and $k\ge\ell\ge1$, it is sufficient to prove the inequality
\begin{equation}\label{U-V-k-com}
\frac{U_{k+1}}{U_{k}}\le\frac{V_{k+1}(1)}{V_k(1)}
\end{equation}
for $k\ge4$, which is equivalent to
\begin{multline*}
\bigl(k^2-3 k-12\bigr)2^{k+1} 3^{k+4}
+\bigl(k^2+10 k+20\bigr)3^{k+4} +k^2+6 k+4\\
+\bigl[\bigl(k^2+13 k+20\bigr)2^{k+5} -\bigl(k^4+16 k^3+118 k^2+435 k+540\bigr)\bigr]2^{k+1}\ge0
\end{multline*}
for $k\ge4$. Since
\begin{multline*}
\bigl(k^2+13 k+20\bigr)2^{k+5} - \bigl(k^4+16 k^3+118 k^2+435 k+540\bigr)\\
>2^{5} \bigl(k^2+13 k+20\bigr)\sum_{\ell=0}^2\binom{k}{\ell}- \bigl(k^4+16 k^3+118 k^2+435 k+540\bigr)\\
=15 k^4+208 k^3+442 k^2+301 k+100
\ge0
\end{multline*}
and $k^2-3 k-12$ is positive for $k\ge6$, the inequality~\eqref{U-V-k-com} is valid for $k\ge6$. By a straightforward computation, it is easy to see that the inequality~\eqref{U-V-k-com} is also valid for $k=4,5$. Therefore, the inequality~\eqref{U-V-k-com} is valid for all $k\ge4$. In conclusion, the inequality~\eqref{theta=k>k-eq} holds for $k\ge4$ and $k\ge\ell\ge1$.
\par
Substituting~\eqref{theta-0-compute} and~\eqref{theta=k>k-eq} into~\eqref{C=k+1=C-k>0-eq} reveals that $C_{k+1}(\beta)-C_k(\beta)>0$ is valid for $k\ge4$ and $0<\beta<1$. Hence, the sequence $C_k(\beta)=\frac{\xi_k(\beta)}{\lambda_k}$ is increasing with respect to $k\ge0$ for $0<\beta<1$. By Lemma~\ref{Ponnusamy-Vuorinen-97-lem}, it follows that the derivative
$$
\frac{\td}{\td u}\biggl[\frac1{G_{\beta}(u)}\biggr] =\frac{\sum_{k=0}^\infty \xi_k(\beta)u^k}{\sum_{k=0}^\infty \lambda_ku^k}
$$
is increasing and that the function $\frac1{G_{\beta}(u)}$ is convex on $(0,\infty)$. The proof of the convexity of the function~\eqref{G-beta-(u)-dfn} is complete.
\par
It is easy to obtain
\begin{equation*}
\lim_{u\to0^+}\frac{\td}{\td u}\biggl[\frac1{G_{\beta}(u)}\biggr] =\lim_{u\to0^+}\frac{\sum_{k=0}^\infty \xi_k(\beta)u^k}{\sum_{k=0}^\infty \lambda_ku^k} =\frac{\xi_0(\beta)}{\lambda_0}=\frac{\beta-1}{3}<0
\end{equation*}
and, by Lemma~\ref{BesselI-Asympt-lem},
\begin{align*}
\frac{\td}{\td u}\biggl[\frac1{G_{\beta}(u)}\biggr]
&=\frac{\beta(e^u-1)^2I_3\bigl(2\sqrt{\beta u}\,\bigr)}{e^{u}(e^u-1-u)(\beta u)^{3/2}}
-\frac{(e^u-1)[e^u (u-2)+u+2] I_2\bigl(2\sqrt{\beta u}\,\bigr)}{\beta ue^{u}(e^u-1-u)^2}\\
&\sim\beta\frac{I_3\bigl(2\sqrt{\beta u}\,\bigr)}{(\beta u)^{3/2}}
-\frac1{\beta e^{u}}{I_2\bigl(2\sqrt{\beta u}\,\bigr)}
\to\infty
\end{align*}
as $u\to\infty$ for $0<\beta<1$. Consequently, from its monotonicity on $(0,\infty)$, the derivative $\frac{\td}{\td u}\bigl[\frac1{G_{\beta}(u)}\bigr]$ has a unique zero, and so the function $\frac1{G_{\beta}(u)}$ has a unique minimum, and so the positive function $G_{\beta}(u)$ has a unique maximum, on $(0,\infty)$.
The proof of the unimodality of the function~\eqref{G-beta-(u)-dfn} is complete.
\end{proof}

\section{Proof of Theorem~\ref{H(t)-degree=4-thm}}

With the help of Theorem~\ref{lower-beseel=5-thm}, we now start off to prove Theorem~\ref{H(t)-degree=4-thm}.
\par
If the function $t^qH_{1,1}(t)$ is completely monotonic on $(0,\infty)$, then its first derivative is non-positive, that is,
\begin{equation*}
-t^{q-2} \bigl\{qt\bigl[\psi'(t)-e^{1/t}+1\bigr]+t^2\psi''(t)+e^{1/t}\bigr\}\le0,
\end{equation*}
which can be formulated as
\begin{equation*}
q\le\frac{t^2\psi''(t)+e^{1/t}}{t\bigl[e^{1/t}-\psi'(t)-1\bigr]}\triangleq p(t).
\end{equation*}
By virtue of~\eqref{asymptotic-polypsi} for $n=1,2$ and the expansion
\begin{equation*}
e^{1/t}=1+\sum_{m=1}^\infty\frac1{m!}\frac1{t^m},\quad t\ne0,
\end{equation*}
we have
\begin{align*}
p(t)&\sim\frac{\sum_{k=0}^3\frac1{k!t^k}+O\bigl(\frac1{t^3}\bigr) -t^2\bigl[\frac1{t^2}+\frac1{t^3}+\frac1{2t^4}-\frac1{6t^6}+O\bigl(\frac1{t^6}\bigr)\bigr]} {t\bigl\{\bigl[\sum_{k=0}^4\frac1{k!t^k}+O\bigl(\frac1{t^4}\bigr)\bigr] -\bigl[\frac1t+\frac1{2t^2}+\frac1{6t^3} -\frac1{30t^5}+O\bigl(\frac1{t^5}\bigr)\bigr]-1\bigr\}}\\
&\sim\frac{\frac1{3!t^3}+O\bigl(\frac1{t^3}\bigr)} {t\bigl[\frac1{4!t^4}+O\bigl(\frac1{t^4}\bigr)\bigr]}
\to4
\end{align*}
as $t\to\infty$. This implies that
\begin{equation}\label{H(t)degree<4}
\cmdeg{t}[H_{1,1}(t)]\le4.
\end{equation}
\par
By the integral representation~\eqref{h(t)-int-repres}, the formula~\eqref{Beseel-der-recur}, the definition of $H_{1,1}(t)$, and integration by part, we have
\begin{align}
t^4H_{1,1}(t)&=t^4\int_0^\infty \biggl[\frac{I_1\bigl(2\sqrt{u}\,\bigr)}{\sqrt{u}\,} -\frac{u}{1-e^{-u}}\biggr] e^{-tu}\td u\notag\\
&=-t^3\int_0^\infty \Biggl[\frac{I_1\bigl(2\sqrt{u}\,\bigr)}{\sqrt{u}\,} -\frac{u}{1-e^{-u}}\Biggr] \frac{\td e^{-tu}}{\td u}\td u\notag\\
&=t^3\int_0^\infty \frac{\td}{\td u} \Biggl[\frac{I_1\bigl(2\sqrt{u}\,\bigr)}{\sqrt{u}\,} -\frac{u}{1-e^{-u}}\Biggr]e^{-tu}\td u\notag\\
&=t^2\int_0^\infty \frac{\td^2}{\td u^2} \Biggl[\frac{I_1\bigl(2\sqrt{u}\,\bigr)}{\sqrt{u}\,} -\frac{u}{1-e^{-u}}\Biggr]e^{-tu}\td u\notag\\
&=t\int_0^\infty \frac{\td^3}{\td u^3} \Biggl[\frac{I_1\bigl(2\sqrt{u}\,\bigr)}{\sqrt{u}\,} -\frac{u}{1-e^{-u}}\Biggr]e^{-tu}\td u\notag\\
&=\frac1{24}+\int_0^\infty \frac{\td^4}{\td u^4} \Biggl[\frac{I_1\bigl(2\sqrt{u}\,\bigr)}{\sqrt{u}\,} -\frac{u}{1-e^{-u}}\Biggr]e^{-tu}\td u\notag\\
&=\frac1{24}+\int_0^\infty \biggl[\frac{I_5\bigl(2\sqrt{u}\,\bigr)}{u^{5/2}} -\biggl(\frac{u}{1-e^{-u}}\biggr)^{(4)}\biggr]e^{-tu}\td u. \label{(3.5)-frac1(24)-eq}
\end{align}
As a result, by Lemma~\ref{Bernstein-Widder-Theorem-12b} and by the inequality~\eqref{conj-bessel-ineq} for $k=5$, we conclude that the function $t^4H_{1,1}(t)$ is completely monotonic on $(0,\infty)$. So, by the definition of completely monotonic degrees, we have
\begin{equation}\label{H(t)degree>4}
\cmdeg{t}[H_{1,1}(t)]\ge4.
\end{equation}
Combining~\eqref{H(t)degree<4} and~\eqref{H(t)degree>4} yields~\eqref{H(t)degree=4}.
\par
As deducing the integral represntation~\eqref{h(t)-int-repres}, we can derive that
\begin{equation}
H_{\alpha,\beta}(z)=\int_0^{\infty}\biggl[\alpha\beta\frac{I_1\bigl(2\sqrt{\beta u}\,\bigr)}{\sqrt{\beta u}\,}-\frac{u}{1-e^{-u}}\biggr]e^{-zu}\td u
\end{equation}
for $\Re z>0$.
Employing~\eqref{Beseel-der-recur} and integrating in part as in~\eqref{(3.5)-frac1(24)-eq} yield
\begin{align*}
tH_{\alpha,\beta}(t)&=-\int_0^{\infty}\biggl[\alpha\beta\frac{I_1\bigl(2\sqrt{\beta u}\,\bigr)}{\sqrt{\beta u}\,}-\frac{u}{1-e^{-u}}\biggr]\frac{\td e^{-tu}}{\td u}\td u\\
&=\alpha\beta-1+\int_0^{\infty}\biggl[\alpha\beta^2\frac{I_2\bigl(2\sqrt{\beta u}\,\bigr)}{\beta u} -\biggl(\frac{u}{1-e^{-u}}\biggr)'\biggr]e^{-tu}\td u.
\end{align*}
Consequently,
\begin{enumerate}
\item
when $\alpha\beta=1$, the function $tH_{\alpha,\beta}(t)$ becomes
\begin{equation*}
tH_{\alpha,\beta}(t)=\int_0^{\infty}\biggl[\beta\frac{I_2\bigl(2\sqrt{\beta u}\,\bigr)}{\beta u} -\biggl(\frac{u}{1-e^{-u}}\biggr)'\biggr]e^{-tu}\td u
\end{equation*}
and, by integration by part and the recursion~\eqref{Beseel-der-recur},
\begin{align*}
t^2H_{\alpha,\beta}(t)&=-\int_0^{\infty} \biggl[\beta\frac{I_2\bigl(2\sqrt{\beta u}\,\bigr)}{\beta u} -\biggl(\frac{u}{1-e^{-u}}\biggr)'\biggr]\frac{\td e^{-tu}}{\td u}\td u\\
&=\frac{\beta-1}2+\int_0^{\infty}\biggl[\beta\frac{I_3\bigl(2\sqrt{\beta u}\,\bigr)}{(\beta u)^{3/2}} -\biggl(\frac{u}{1-e^{-u}}\biggr)''\biggr]e^{-tu}\td u;
\end{align*}
\begin{enumerate}
\item
if $\beta>1$, by virtue of the fact that the function $\frac{I_3(2u)}{u^3}$ is strictly increasing on $(0,\infty)$ and by the inequality~\eqref{conj-bessel-ineq} for $k=3$, it is not difficult to see that the completely monotonic degree of the function $H_{1/\beta,\beta}(t)$ for $\beta>1$ is $2$;
\item
if $0<\beta<1$, by the necessary condition~\eqref{beta<1-condition}, the function $H_{\alpha,\beta}(t)$ is not completely monotonic;
\item
if $\beta=1$, the discussing question goes back to the proof of~\eqref{H(t)degree=4};
\end{enumerate}
\item
when $\alpha\beta>1$ and
\begin{equation}
\alpha\beta^2\ge\frac{\beta u}{I_2\bigl(2\sqrt{\beta u}\,\bigr)}\biggl(\frac{u}{1-e^{-u}}\biggr)'=G_{\beta}(u),
\end{equation}
the completely monotonic degree of $H_{\alpha,\beta}(t)$ is $1$; By virtue of the monotonicity and unimodality of the function~\eqref{G-beta-(u)-dfn} obtained in Theorem~\ref{lower-beseel=5-thm}, the quantity~\eqref{H-alpha-beta(t)degree=1} follows.
\end{enumerate}
The proof of Theorem~\ref{H(t)-degree=4-thm} is complete.

\section{Remarks}

Finally we list some remarks on something to do with our lemmas and theorems.

\begin{rem}
We note that Lemma~\ref{Ponnusamy-Vuorinen-97-lem} has been generalized in~\cite[Lemma~2.2]{Koumandos-Pedersen-287}.
\end{rem}

\begin{rem}
The function $F_3(u)$ defined by~\eqref{F=3(u)-dfn} can also be decomposed as
\begin{equation*}
F_3(u)=f_1(u)+f_2(u)+f_3(u),
\end{equation*}
where
\begin{gather*}
f_1(u)=[10u(e^u-261)+3966]e^{2u},\quad
f_2(u)=\bigl(69e^{2u}-7119u-4035\bigr)e^u,\\
f_3(u)=3249e^{2u}-793u-3249,
\end{gather*}
and these three functions are all increasing respectively on the intervals $[5,\infty)$, $[3,\infty)$, and $[0,\infty)$. Then it follows that $F_3(u)$ is increasing and positive on $[5,\infty)$.
\end{rem}

\begin{rem}
In the draft of this manuscript, we ever used Theorem~2 in~\cite[p.~22]{exp-beograd} to prove the positivity of the function $F_3(u)$ defined in~\eqref{F=3(u)-dfn} on $(0,6)$. But, the inequality~(14) stated in~\cite[p.~22, Theorem~2]{exp-beograd}, and then the inequality~(18) in~\cite[p.~22]{exp-beograd}, is wrong. So we have to give up using~\cite[p.~22, Theorem~2]{exp-beograd} to prove the inequality~\eqref{conj-bessel-ineq}.
\par
By the way, we can reformulate~\cite[Theorem~1]{exp-beograd} as follows. For $x\in[0,b]$ and $n\ge0$, we have
\begin{equation}\label{ineq-14-15-method-exp-comb}
0\ge e^x -S_n(x)-\alpha_n(b)x^{n+1} \ge\frac{(n+1)!\alpha_n(b)-e^b}{(n+1)!(n+1)b}(b-x)x^{n+1},
\end{equation}
where
\begin{equation}
S_n(x)=\sum_{k=0}^n\frac{x^k}{k!}\quad \text{and}\quad \alpha_n(b)=\frac{e^b-S_n(b)}{b^{n+1}}.
\end{equation}
The equalities in~\eqref{ineq-14-15-method-exp-comb} are valid if and only if $x=0,b$. This theorem was proved once again in~\cite{Jin-Zhang-Huzhou-09} and was collected in the monograph~\cite[p.~290]{kuang-3rd} and its older and subsequently revised version.
\end{rem}

\begin{rem}
By Descartes' Sign Rule, it follows that
\begin{enumerate}
\item
the polynomial $P_1(m)=m^4+36m^3+323m^2+12m-756$ has one possible positive zero; since $P_1(0)=-756$ and $P_1(2)=864$, this zero belongs to the interval $(0,2)$, so $P_1(m)>0$ for $m\ge2$;
\item
the polynomial $P_2(m)=m^4+10m^3-99m^2+24m+252$ has two possible positive zeros; since $P_2(0)=252$, $P_2(3)=-216$, and $P_2(6)=288$, these two zeros locate in the interval $(0,6)$, so $P_2(m)>0$ for $m\ge6$.
\end{enumerate}
Furthermore, we have
\begin{align*}
&\quad m(m+1)9^{m+6}-2\bigl(4m^4+52m^3+207m^2+327m+288\bigr)3^{m+5}\\
&=\bigl[m(m+1)3^{m+7}-2\bigl(4m^4+52m^3+207m^2+327m+288\bigr)\bigr]3^{m+5}\\
&>\Biggl[m(m+1)\sum_{\ell=0}^3\binom{m+7}{\ell}2^\ell -2\bigl(4m^4+52m^3+207m^2+327m+288\bigr)\Biggr]3^{m+5}\\
&=\bigl(4 m^5+58 m^4+278 m^3+407 m^2-825 m-1728\bigr)3^{m+4}
\end{align*}
and, by Descartes' Sign Rule, the polynomial
\begin{equation*}
P_3(m)=4m^5+58m^4+278m^3+407m^2-825m-1728
\end{equation*}
has one possible positive zero. Since $P_3(0)=-1728$ and $P_3(1)=1530$, it follows that $P_3(m)$ is positive for $m\ge2$. Consequently, we obtain that $\mathcal{C}(m)$ defined in~\eqref{M-m(k)-dfn} is positive for $m\ge6$. Considering that
\begin{align*}
\mathcal{C}(0)&=181440,& \mathcal{C}(1)&=10160640, & \mathcal{C}(2)&=252316512,\\
\mathcal{C}(3)&=4549288320,& \mathcal{C}(4)&=68981774400,& \mathcal{C}(5)&=939390217920,
\end{align*}
we conclude that $\mathcal{C}(m)$ are positive for all nonnegative integers $m\ge0$.
\par
By similar argument to above, we can determine that $\mathcal{B}(m)<0$ for all nonnegative integers $m\ge0$.
\end{rem}

\begin{rem}
In order to prove
\begin{equation}\label{u+6-exp-4-deriv}
\frac{u+6}{720}\ge\biggl(\frac{u}{1-e^{-u}}\biggr)^{(4)},
\end{equation}
we write
\begin{equation*}
\frac{u}{1-e^{-u}}=u\sum_{k=0}^\infty e^{-ku}
\end{equation*}
from which we obtain
\begin{equation*}
\biggl(\frac{u}{1-e^{-u}}\biggr)^{(4)}=\sum_{k=1}^\infty k^3(ku-4)e^{-ku}.
\end{equation*}
So for $u\ge a\ge4$
\begin{equation*}
\biggl(\frac{u}{1-e^{-u}}\biggr)^{(4)}\le\sum_{k=1}^\infty k^3(ku-4)e^{-ka}=uK_4(a)-4K_3(a),
\end{equation*}
where
\begin{equation*}
K_\ell(a)=\sum_{k=1}^\infty k^\ell e^{-ka}.
\end{equation*}
These functions can be calculated for small values of $\ell$ and $K_4(7)=0.0009\dotsm<\frac1{720}$. Therefore, the inequality~\eqref{u+6-exp-4-deriv} holds for $u\ge7$. As a result, it is sufficient to prove the inequality~\eqref{u+6-exp-4-deriv} on the interval $[0,7]$.
\end{rem}

\begin{rem}
By a result in~\cite{Temme-96-book} (or see~\cite[p.~35, (3)]{Koumandos-Pedersen-09-JMAA}), we have
\begin{equation}\label{V=n(u)-dfn}
\frac{u}{1-e^{-u}}=1+\frac{u}2+\sum_{k=1}^n \frac{B_{2k}}{(2k)!}u^{2k}+(-1)^nu^{2(n+1)}V_n(u),
\end{equation}
where
\begin{equation}
V_n(u)=\sum_{k=1}^\infty\frac2{(u^2+4\pi^2k^2)(2\pi k)^{2n}}
\end{equation}
and it was proved in~\cite[Lemma~2.3]{Koumandos-Pedersen-09-JMAA} that
\begin{equation}\label{V=n(u)>0-eq}
\bigl[u^{2\ell}V_n(u)\bigr]^{(k)}\ge0
\end{equation}
for $u>0$ and $0\le k\le\ell$. By~\eqref{V=n(u)-dfn} for $n=1$, it follows that
\begin{equation}
u^4V_1(u)=1+\frac{u}2+\frac{u^2}{12}-\frac{u}{1-e^{-u}},
\end{equation}
hence, by~\eqref{polygamma} for $n=1$,
\begin{equation}
\int_0^\infty u^4V_1(u)e^{-ux}\td u=\frac1x+\frac1{2x^2}+\frac1{6x^3}-\psi'(x).
\end{equation}
Since the inequality~\eqref{V=n(u)>0-eq} holds for $0\le k\le2$, we obtain by~\cite[Theorem~1.3]{Koumandos-Pedersen-09-JMAA} that
\begin{equation*}
x^2\biggl[\frac1x+\frac1{2x^2}+\frac1{6x^3}-\psi'(x)\biggr]
\end{equation*}
is completely monotonic on $(0,\infty)$. If we add the completely monotonic function
\begin{equation*}
x^2\biggl(e^{1/x}-1-\frac1x-\frac1{2x^2}-\frac1{6x^3}\biggr)
\end{equation*}
to the above, we find that the function
\begin{equation*}
x^2\bigl[e^{1/x}-1-\psi'(x)\bigr]
\end{equation*}
is completely monotonic on $(0,\infty)$, that is,
\begin{equation}
\cmdeg{t}[H_{1,1}(t)]\ge2.
\end{equation}
\end{rem}

\begin{rem}
Now, if we use~\eqref{V=n(u)-dfn} for $n=3$, we obtain similarly that
\begin{equation}
\int_0^\infty u^8V_3(u)e^{-ux}\td u=\frac1x+\frac1{2x^2}+\frac1{6x^3} -\frac1{30x^5}+\frac1{42x^7}-\psi'(x).
\end{equation}
Since the inequality~\eqref{V=n(u)>0-eq} holds for $0\le k\le4$, we acquire by~\cite[Theorem~1.3]{Koumandos-Pedersen-09-JMAA} that
\begin{equation*}
x^4\biggl[\frac1x+\frac1{2x^2}+\frac1{6x^3} -\frac1{30x^5}+\frac1{42x^7}-\psi'(x)\biggr]
\end{equation*}
is completely monotonic on $(0,\infty)$. If we add the completely monotonic function
\begin{equation*}
x^4\biggl(e^{1/x}-1-\frac1x-\frac1{2x^2}-\frac1{6x^3}-\frac1{4!x^4} -\frac1{5!x^5}-\frac1{6!x^6}-\frac1{7!x^7}\biggr)
\end{equation*}
to the above, we gain that the function
\begin{equation*}
x^4\bigl[e^{1/x}-1-\psi'(x)\bigr]-\frac1{24}-\frac1{24x}-\frac1{6!x^2}+\frac{17}{6!x^3}
\end{equation*}
is completely monotonic on $(0,\infty)$. By further adding three completely monotonic terms we finally earn that the function
\begin{equation*}
x^4\bigl[e^{1/x}-1-\psi'(x)\bigr]-\frac1{24}+\frac{17}{6!x^3}
\end{equation*}
is completely monotonic on $(0,\infty)$.
\par
If one can manage to remove the last term $\frac{17}{6!x^3}$, then the first result~\eqref{H(t)degree=4} in Theorem~\ref{H(t)-degree=4-thm} follows.
\end{rem}

\begin{rem}
Motivated by properties of the functions $F_{\beta}(u)$ and $G_{\beta}(u)$ defined in~\eqref{F-beta-(u)-dfn} and~\eqref{G-beta-(u)-dfn} respectively, we conjecture that
\begin{enumerate}
\item
when $3\le k\le 5$ and $\beta\ge1$, the function
\begin{equation}\label{conj-bessel-ineq-k}
H_{k,\beta}(u)=\frac{(\beta u)^{k/2}}{I_k\bigl(2\sqrt{\beta u}\,\bigr)}\biggl(\frac{u}{1-e^{-u}}\biggr)^{(k-1)}
\end{equation}
is decreasing on $(0,\infty)$;
\item
when $3\le k\le 5$ and $0<\beta<1$, the function $H_{k,\beta}(u)$ is unimodal on $(0,\infty)$;
\item
when $3\le k\le 5$ and $0<\beta<1$, the function $\frac1{H_{k,\beta}(u)}$ is convex on $(0,\infty)$.
\end{enumerate}
\end{rem}

\begin{rem}
We conjecture that for all $k\ge6$ the inequality~\eqref{conj-bessel-ineq} does not hold on $(0,\infty)$.
\end{rem}

\end{document}